\documentclass{article}
\usepackage{amsmath}
\usepackage{latexsym}
\usepackage{amssymb}

\newtheorem{theorem}{Theorem}

\newtheorem{lemma}[theorem]{Lemma}

\newtheorem{remark}[theorem]{Remark}

\title{On Decay of Solutions to Nonlinear Schr\"odinger Equations}

\author{Alexander Pankov\\
Mathematics Department\\
College of William and Mary\\
Williamsburg, VA 23187--8795\\
e-mail: {\tt pankov@member.ams.org}}
\date{}

\begin{document}

\maketitle

\begin{abstract} We present general results on exponential decay of finite energy solutions to
stationary nonlinear Schr\"odinger equations.

\vspace{2ex}

\end{abstract}

In this note we consider the equation
\begin{equation}\label{e1}
-\Delta u +V(x)u=f(x,u), \quad x\in \mathbb{R}^n\,
\end{equation}
and, under rather general assumptions, derive exponential decay estimates for its solutions.

We suppose that

\vspace{1ex}

$(i)$ {\em The potential\ \  $V$ belongs to $L^{\infty}_{\mathrm{loc}}(\mathbb{R}^n)$ and is bounded
below, i.e. $V(x)\geq -c_0$ for some $c_0\in \mathbb{R}$.}

\vspace{1ex}

\noindent Under assumption $(i)$ the left hand side of equation (\ref{e1}) defines a self-adjoint
operator in $L^2(\mathbb{R}^n)$ denoted by $H$. The operator $H$ is bounded below. We suppose that

\vspace{1ex}

$(ii)$ {\em The essential spectrum $\sigma_{\mathrm{ess}}(H)$ of the operator $H$ does not contain the
point $0$.}

\vspace{1ex}

\noindent Note, however, that $0$ can be an eigenvalue of finite multiplicity.

The nonlinearity $f$ is supposed to satisfy the following assumption.

\vspace{1ex}

$(iii)${\em The function $f(x, u)$ is a Carath\'eodory function, i.e. it is Lebesgue measurable with
respect to $x\in\mathbb{R}^n$ for all $u\in\mathbb{R}$ and continuous with respect to $u\in\mathbb{R}$
for almost all $x\in\mathbb{R}^n$. Furthermore,
\begin{equation}\label{e2}
|f(x,u)|\leq c(1+|u|^{p-1})\,,\quad x\in\mathbb{R}^n\, u\in\mathbb{R}\,,
\end{equation}
with $c>0$ and $2\leq p<2^*$, where
$$
2^*=\begin{cases}\displaystyle\frac{2n}{n-2}\quad &\text{if } n\geq 3\,,\\ \infty\quad &\text{if } n=1,
2\,,
\end{cases}
$$
and
$$
\lim_{u\to 0}\mathrm{ess \, sup}_{x\in\mathbb{R}^n}\frac{|f(x,u)|}{|u|} =0\,.
$$}

\vspace{1ex}

Let $E$ denote the form domain of the operator $H$. It is well-known that
$$
E=\{u\in H^1(\mathbb{R}^n)\,:\, (V(x)+c_0+1)u(x)\in L^2(\mathbb{R}^n)\}
$$
where $c_0$ is the constant from assumption $(i)$. Moreover, $E$ carries a natural Hilbert space
structure and is continuously embedded into the space $H^1(\mathbb{R}^n)$. In what follows we consider
only weak solutions that belong to the space $E$. A function $u\in E$ is a weak solution of
equation~(\ref{e1}) if for all $v\in E$ the following integral identity
\begin{equation}\label{e3}
\int_{\mathbb{R}^n}(\nabla u(x)\cdot\nabla v(x)+V(x)u(x)v(x)-f(x,u(x))v(x))\,dx
\end{equation}
is satisfied. Actually, it is sufficient to check identity~(\ref{e3}) only for functions $v$ that belong
to the space $C_0^{\infty}(\mathbb{R}^n)$ of all finitely supported infinitely differentiable functions.
Also we note that, due to the Sobolev embedding theorem, the term $f(x, u(x))v(x)$ in (\ref{e3}) is
integrable for all $u,v\in H^1(\mathbb{R}^n)$.

\vspace{2ex}

There is a number of results on exponential decay of solutions to equation (\ref{e1}) (see \cite{bpw,
b-l, fu-oz, pank05, pank06, stu98}). Most of them, except \cite{pank05, pank06}, deal with the case when
$0$ is below the essential spectrum of $H$. However, the case when $0$ is in a spectral gap is extremely
important for applications \cite{pank05}. Here we present a rather general result in this direction.

We exploit the following rather simple idea. Suppose that $u\in E$ is a nontrivial, {\em i. e.} $u\neq
0$, solution to equation (\ref{e1}). Set
$$
W(x)=\begin{cases}\displaystyle\frac{f(x,u(x))}{u(x)}\quad &\text{if } u(x)\neq 0\,,\\
0 \quad &\text{if } u(x)=0\,,
\end{cases}
$$
Then equation~(\ref{e1}) can be represented as
\begin{equation}\label{e4}
(H+W)u=0\,.
\end{equation}
This means that $u$ is an eigenfunction of the operator $H+W$ with zero eigenvalue. Now if the
multiplication operator by $W$ is relatively compact with respect to the operator $H$, then
$$
\sigma_{\mathrm{ess}}(H+W)=\sigma_{\mathrm{ess}}(H+W)\,.
$$
Hence, 0 is an eigenvalue of $H+W$ of finite multiplicity and $u$ a corresponding eigenfunction (see,
e.g. \cite{h-s, r-s}). Now an exponential decay of $u$ can be read off from any well-known result about
eigenfunctions. More precisely, the solution $u$ has exactly the same decay as an eigenfunction that
corresponds to an eigenvalue introduced into a spectral gap by a decaying perturbation of the potential
$V$.

Thus, the only we need is to verify that the multiplication operator by the function $W$ is relatively
compact with respect to $H$. The key point is the following
\begin{lemma}\label{l1}
Under assumptions $(i)$--$(iii)$ suppose that $u\in E$ is a solution of equation (\ref{e1}). Then $u$ is
a continuous function and
$$
\lim_{x\to\infty}u(x)=0\,.
$$
\end{lemma}

We postpone the proof of the lemma and first present main results.

\vspace{2ex}

Due to a well-known result (see,  e.g., Theorem~8.3.1 of \cite{pank06}), Lemma~\ref{l1} implies that the
multiplication operator by $W$ is a relatively compact perturbation of the operator $H$. Making use of
Theorem~C.3.4, \cite{si82}, we obtaine

\begin{theorem}\label{t1} Assume $(i)$--$(iii)$.
Let $u\in E$ be a solution of equation (\ref{e1}). Then there exists $\alpha_0>0$ such that for every
$\alpha <\alpha_0$ we have
\begin{equation}\label{e5}
|u(x)|\leq C \exp(-\alpha |x|)\,
\end{equation}
with some $C=C_{\alpha}>0$.
\end{theorem}

\begin{remark}\label{r0} An interesting case is when the potential is periodic. If, in addition, the
nonlinearity is superlinear, i.e. $f(\cdot,u)\geq c|u|^{p-1}$ at infinity, with $p>2$, the result of
Theorem~\ref{t1} is announced in \cite{pank05}. However, assumption $(iii)$ allows asymptotically linear
nonlinearities ($p=2$). As consequence, solutions found in \cite{li-sz} decay exponentially fast.
\end{remark}

\begin{remark}\label{r1}
The value of $\alpha_0$ can be estimated in terms of the distance between $0$ and
$\sigma_{\mathrm{ess}}(H)$ (see, e.g., \cite{hisl}).
\end{remark}

Now we consider the case when $\sigma_{\mathrm{ess}}(H)=\varnothing$, {\em i.e.} the spectrum of $H$ is
discrete. This is so if,  e.g,
$$
\lim_{|x|\to\infty}V(x)=\infty\,.
$$
For a necessary and sufficient condition for the discreteness of spectrum see \cite{Koshu99} and
references therein. In this case Theorem~C.3.3 of \cite{si82} implies

\begin{theorem}\label{t2} Under assumptions $(i)$--$((iii)$, suppose that the spectrum of $H$ is
discrete. Let $u\in E$ be a solution of equation (\ref{e1}). Then fore every $\alpha>0$ there exists
$C=C_{\alpha}>0$ such that
\begin{equation}\label{e6}
|u(x)|\leq C \exp(-\alpha |x|)\,.
\end{equation}
\end{theorem}

Having an additional information about the behavior of $V$ at infinity one can refine the result of
Theorem~\ref{t2}. For instance, making use of Theorem~3.3, \cite{b-s}, we obtain

\begin{theorem}\label{t3}
In addition to assumptions $(i)$--$(iii)$, suppose that
\begin{equation}\label{e7}
V(x)\geq \gamma |x|^{\beta}-\gamma_0\,
\end{equation}
with $\gamma>0$, $\gamma_0\geq 0$ and $\beta>0$. Then for any solution $u\in E$ of equation \ref{e1} we
have that
\begin{equation}\label{e8}
|u(x)\leq C \exp(-a|x|^{\frac{\beta}{2}+1})\,,
\end{equation}
with some $C>0$ and $a>0$.
\end{theorem}

Now we prove Lemma~\ref{l1}.

{\em Proof of Lemma~\ref{l1}\,}. Case $n=1$ is trivial because any function from $H^1(\mathbb{R})$ is
continuous and vanishes at infinity.

Now we consider case $n\geq 3$. We use a version of the well-known bootstrap argument as follows.

Equation (\ref{e1}) can be rewritten as
   \begin{equation}\label{e9}
      (-\Delta +V+c_0+1)\,u=(c_0+1)\,u+f(x,u)\,.
   \end{equation}
   Note that by the Sobolev embedding   $u\in L^{2^*}(\mathbb{R}^n)$.

   Suppose now that $u\in L^\infty(\mathbb{R}^n)+L^r(\mathbb{R}^n)$ with $r\ge 2^*$. Then the second term in
   the right hand side of (\ref{e9}) belongs to $L^s(\mathbb{R}^n)$ with $s=r/(p-1)$, while the first
   term belongs to $L^r(\mathbb{R}^n)$. Obviously, $r>s$.

   Let
      $$A=\Big\{x\in\mathbb{R}^n\;:\; \big|u(x)\big|\ge 1\Big\}$$
   and $B=\mathbb{R}^n\setminus A$. It is easy that ${\rm meas}\,(A)<\infty$.  Denote by $\chi_A$ and $\chi_B$ the characteristic functions of the sets
   $A$ and $B$ respectively, i.~e. $\chi_A=1$ on $A$, $\chi_A=0$ on $B$ and $\chi_B=1-\chi_A$.

   Let
   \begin{gather*}
      H_1=-\Delta+V(x)+c_0+1\,,\\
      h_0(x)=\chi_B(x)\,\Big[(c_0+1)\,u(x)+f(x,u)\Big]\,,
   \end{gather*}
   and
      $$  h_1(x)=\chi_A(x)\,\Big[(c_0+1)\,u(x)+f(x,u)\Big]\,.$$
   Equation (\ref{e9}) becomes
   \begin{equation}\label{e10}
      H_1u=h_0(x)+h_1(x).
   \end{equation}
Obviously, $h_0\in L^\infty(\mathbb{R}^n)$, while assumption $(iii)$ implies that $h_1\in
L^s(\mathbb{R}^n)$.

   The operator $H_1$ is positive definite and satisfies the assumptions of the Sobolev
   estimate theorem for Schr\"odinger operators (see Theorem~B.2.1 of \cite{si82}). Hence, we obtain from
   (\ref{e10}) that $u=u_0+u_1$, where
      $$ u_0=H_1^{-1}h_0\in L^\infty(\mathbb{R}^n)$$
    and
      $$ u_1=H_1^{-1}h_1\in L^q(\mathbb{R}^n)$$
   for every $q$ such that
   \begin{equation}\label{e11}
      \frac 1s-\frac 1q<\frac 2n\,.
   \end{equation}
   If $s>n/2$, {\em i.~e.}
   \begin{equation}\label{e12}
       r>\frac{n(p-1)}2\,,
   \end{equation}
   we can take $q=\infty$ to obtain  $u_1\in L^\infty(\mathbb{R}^n)$ and, hence, $u=u_0+u_1\in
   L^\infty(\mathbb{R}^n)$.

   Otherwise, take
      $$ q=\frac r{1-\delta}\,,$$
   where
      $$\delta=\frac 4{n-2}-(p-2)-\varepsilon$$
   and $\varepsilon>0$ is arbitrary small. Note that
      $$ p-2<\frac 4{n-2}$$
   because $p<2^*$. Since $r\ge 2^*$, we have that
   \begin{align*}
      \frac 1s-\frac 1q &= \frac{p-1}r-\frac{1-\delta}r=\frac 1r\left(\frac 4{n-2}-\varepsilon\right)
                         \le \frac 1{2^*}\left(\frac 4{n-2}-\varepsilon\right)\\
                        &=\frac 2n-\frac{n-2}{2n}\,\varepsilon\,.
   \end{align*}
   Therefore, $q$ satisfies (\ref{e12}) and $u_1\in L^\infty(\mathbb{R}^n)$. Hence,
      $$u=u_0+u_1\in L^\infty(\mathbb{R}^n)+L^q(\mathbb{R}^n)\,.$$
It is not difficult to verify  that $\chi_Au\in L^q(\mathbb{R}^n)$ and $h_1\in
L^{q/(p-1)}(\mathbb{R}^n)$.

\vspace{3ex}

   Now starting with $r=r_0=2^*$, we can iterate the previous procedure. Thus, we have that $u\in
   L^\infty(\mathbb{R}^n)+L^{r_k}(\mathbb{R}^n)$, where
      $$ r_k=\frac{2^*}{(1-\delta)^k}\,.$$
   Let $k$ be so large that
      $$ r_k>\frac{n(p-1)}2$$
   (see (\ref{e12})). Then
      $$ s_k=\frac{r_k}{p-1}>\frac n2$$
   and we can apply the Sobolev estimate of Theorem~B.2.1 \cite{si82}
    with $s=s_k$ and $q=\infty$. Therefore, $u_1\in
   L^\infty(\mathbb{R}^n)$ and, hence, $u=u_0+u_1\in L^\infty(\mathbb{R}^n)$.

   This implies immediately that the additional potential $W$ belongs to
   $L^{\infty}_{\mathrm{loc}}(\mathbb{R}^n)$ and is bounded below. Therefore, due to Theorems~C.1.1 and
   C.3.1 of \cite{si82} the result follows.

   Case $n=2$ is simpler. By the Sobolev embedding, $u\in L^r(\mathbb{R}^n)$ for arbitrarily large $r$.
   Hence, the previous argument shows that $u\in L^\infty(\mathbb{R}^n)$ and we are done. \hfil$\Box$

\end{document}